\documentclass[a4paper,12pt,reqno]{amsart}

\usepackage{latexsym}
\usepackage{amsfonts,amsmath,amssymb}
\usepackage{euscript}
\usepackage{tikz}

\newcommand{\gauss}[2]{\genfrac{[}{]}{0pt}{}{#1}{#2}}

\allowdisplaybreaks

\theoremstyle{remark}

\begin{document}
\title{Factorizations related to the reciprocal Pascal matrix}

\author[Helmut Prodinger]{Helmut Prodinger}
\address{Department of Mathematics, University of Stellenbosch 7602,
Stellenbosch, South Africa}
\email{hproding@sun.ac.za}
\thanks{The author was supported by an incentive grant of the National Research Foundation of South Africa.}

\begin{abstract}The reciprocal Pascal matrix has entries $\binom{i+j}{j}^{-1}$. Explicit formul\ae{} for its LU-decomposition, the LU-decomposition of its inverse, and some  related matrices are obtained. For all results, $q$-analogues are also presented.
\end{abstract}

\maketitle

\section{Introduction}

Recently, there has been some interest in the \emph{reciprocal Pascal matrix} $M$, defined by
\begin{equation*}
M_{i,j}=\binom{i+j}{j}^{-1};
\end{equation*}
the indices start here for convenience with $0,0$, and the matrix is either infinite or has
$N$ rows and columns, depending on the context.

Richardson~\cite{Richardson14} has provided the decomposition
$S=GMG$, where the diagonal matrix $G$ has entries $G_{i,i}=\binom{2i}{i}$, and $S$ is the \emph{super Catalan matrix}~\cite{Gessel92, HeLiMa14} with entries $$S_{i,j}=\frac{(2i)!(2j)!}{i!j!(i+j)!}.$$

We want to give an alternative decomposition of $M$, provided by the LU-decomposition. We will give explicit expressions for $L$ and $U$, defined by $LU=M$, as well as for  $L^{-1}$ and $U^{-1}$.

Since there is also interest in $M^{-1}$, in particular in the integrality of its coefficients, we also provide the LU-decomposition $AB=M^{-1}$, and give expressions for $A$, $B$, $A^{-1}$ and $B^{-1}$.

In the last section, we provide $q$-analogues of these results.

\section{Identities}

The LU-decomposition $M=LU$ is given by
\begin{equation*}
L_{i,j}=\frac{i!i!(2j)!}{(i+j)!(i-j)!j!j!}
\end{equation*}
and
\begin{equation*}
U_{i,j}=\frac{(-1)^ij!j!i!(i-1)!}{(j+i)!(j-i)!(2i-1)!}\qquad \text{for}\ i\ge1.
\end{equation*}
For $i=0$, the formula is $U_{0,j}=1$.

The formula that needs to be proved is
\begin{equation*}
\sum_{0\le k\le \min\{i,j\}}L_{i,k}U_{k,j}=\binom{i+j}{j}^{-1},
\end{equation*}
which is equivalent to 
\begin{equation*}
1+\frac{2i!i!j!j!}{(2i)!(2j)!}\sum_{1\le k\le \min\{i,j\}}
(-1)^k\binom{2i}{i+k}\binom{2j}{j+k}=\binom{i+j}{j}^{-1}.
\end{equation*}
The von Szily identity \cite{Szily94, Gessel92, GeLaFr03} is
\begin{equation*}
\frac{(2i)!(2j)!}{i!j!(i+j)!}=\sum_{k\in\mathbb{Z}}(-1)^k\binom{2i}{i+k}\binom{2j}{j+k},
\end{equation*}
and an equivalent form is, by symmetry,
\begin{equation*}
\frac{(2i)!(2j)!}{i!j!(i+j)!}=\binom{2i}{i}\binom{2j}{j}+2\sum_{k\ge1}(-1)^k\binom{2i}{i+k}\binom{2j}{j+k}.
\end{equation*}
Thus, the identity to be proven is now
\begin{equation*}
\binom{i+j}{j}+\frac{i!j!(i+j)!}{(2i)!(2j)!}\bigg[\frac{(2i)!(2j)!}{i!j!(i+j)!}-\binom{2i}{i}\binom{2j}{j}\biggr]=1,
\end{equation*}
which is obviously correct.

The formula for $L^{-1}$ is for $i\ge j\ge0$:
\begin{equation*}
L^{-1}_{i,j}=\frac{(-1)^{i-j}i!i!(i+j-1)!}{(2i-1)!(i-j)!j!j!}.
\end{equation*}
If necessary ($i=j=0$), this must be interpreted as a limit.

To check this, we consider
\begin{align*}
\sum_{k}&\frac{i!i!(2k)!}{(i+k)!(i-k)!k!k!}
\frac{(-1)^{k-j}k!k!(k+j-1)!}{(2k-1)!(k-j)!j!j!}\\
&=\frac{2i!i!(-1)^{j}}{j!j!}\sum_{j\le k\le i}\frac{k}{(i+k)!(i-k)!}
\frac{(-1)^{k}(k+j-1)!}{(k-j)!}.
\end{align*}
The sum can be evaluated by computer algebra (or otherwise), and the result is indeed $[\![i=j]\!]$, as desired.

The formula for $U^{-1}$ is for  $j\ge i\ge1$
\begin{equation*}
U^{-1}_{i,j}=\frac{(-1)^i(j+i)!(2j)!}{(j-i)!j!(j+i)!(j-1)!i!i!}
\end{equation*}
and for $i=0$:
\begin{equation*}
U^{-1}_{0,j}=	\frac{(2j)!}{j!j!}	.
\end{equation*}

The fact that $\sum_k U_{i,k}U^{-1}_{k,j}=[\![i=j]\!]$ can also be done by computer algebra. Since there are a few cases to be distinguished, it is omitted here.

The LU-decomposition $AB=M^{-1}$ depends on the dimension $N$ and is given by
\begin{equation*}
	A_{i,j}=\frac{(-1)^{i-j}(N-j-1)!j!(N+i-1)!}{i!(N-i-1)!(N+j-1)!(i-j)!},
	\end{equation*}
	\begin{equation*}
	B_{i,j}=\frac{(-1)^{j+N-1}(N+j-1)!}{j!(j-i)!(N-j-1)!i!}.
	\end{equation*}

Since $M^{-1}$ does not have ``nice'' entries, we rather provide formul\ae{} for $A^{-1}$ and 
$B^{-1}$ and prove the identity $B^{-1}A^{-1}=M$ instead. The results are:
\begin{equation*}
A_{i,j}^{-1}=\frac{(N-j-1)!j!(N+i-1)!}{i!(N-i-1)!(N+j-1)!(i-j)!},
\end{equation*}
\begin{equation*}
B_{i,j}^{-1}=\frac{(-1)^{j+N-1}(N-1-i)!j!i!}{(j-i)!(N+i-1)!}.
\end{equation*}
First we prove that these are indeed the inverses. We consider
\begin{align*}
\sum_k&\frac{(-1)^{i-k}(N-k-1)!k!(N+i-1)!}{i!(N-i-1)!(N+k-1)!(i-k)!}
\frac{(N-j-1)!j!(N+k-1)!}{k!(N-k-1)!(N+j-1)!(k-j)!}\\
&=(-1)^{i}\frac{(N+i-1)!(N-j-1)!j!}{(N-i-1)!(N+j-1)!i!}\sum_{j\le k\le i}\frac{(-1)^{k}}{(i-k)!(k-j)!}\\
&=\frac{(N+i-1)!(N-j-1)!j!}{(N-i-1)!(N+j-1)!i!(i-j)!}\sum_{j\le k\le i}(-1)^{i-k}\binom{i-j}{i-k}\\
&=\frac{(N+i-1)!(N-j-1)!j!}{(N-i-1)!(N+j-1)!i!(i-j)!}[\![i=j]\!]=[\![i=j]\!],
\end{align*}
which proves $AA^{-1}=I$. Similarly
\begin{align*}
\sum_k&\frac{(-1)^{k+N-1}(N-1-i)!k!i!}{(k-i)!(N+i-1)!}
\frac{(-1)^{j+N-1}(N+j-1)!}{j!(j-k)!(N-j-1)!k!}\\
&=(-1)^{j}\frac{(N-1-i)!i!(N+j-1)!}{(N+i-1)!j!(N-j-1)!}\sum_k\frac{(-1)^k}{(k-i)!(j-k)!}\\
&=\frac{(N-1-i)!i!(N+j-1)!}{(N+i-1)!j!(N-j-1)!(j-i)!}\sum_k(-1)^{j-k}\binom{j-i}{j-k}\\
&=\frac{(N-1-i)!i!(N+j-1)!}{(N+i-1)!j!(N-j-1)!(j-i)!}[\![i=j]\!]=[\![i=j]\!],
\end{align*}
which proves $B^{-1}B=I$.

Now we compute an entry in $B^{-1}A^{-1}$:
\begin{align*}
\sum_k&\frac{(-1)^{k+N-1}(N-1-i)!k!i!}{(k-i)!(N+i-1)!}
\frac{(N-j-1)!j!(N+k-1)!}{k!(N-k-1)!(N+j-1)!(k-j)!}\\
&=(-1)^{N-1}\frac{(N-1-i)!i!(N-j-1)!j!}{(N+i-1)!(N+j-1)!}\sum_k
\frac{(-1)^k(N+k-1)!}{(k-i)!(N-k-1)!(k-j)!}\\
&=(-1)^{N-1}\frac{i!j!(N-j-1)!}{(N+i-1)!}\sum_k
(-1)^k\binom{N-1-i}{N-1-k}\binom{N+k-1}{N-1+j}\\
&=\frac{i!j!(N-j-1)!}{(N+i-1)!}\sum_k
\binom{i-1-k}{N-1-k}\binom{N+k-1}{N-1+j}\\
&=\frac{i!j!(N-j-1)!}{(N+i-1)!}\sum_k
\binom{i-1-k}{i-N}\binom{N+k-1}{N-1+j}\\
&=\frac{i!j!(N-j-1)!}{(N+i-1)!}
\binom{i-1+N}{i+j}\\
&=\frac{i!j!}{(i+j)!}=M_{i,j},
\end{align*}
as claimed.

Now we use the form $M^{-1}=AB$ and write the $(i,j)$ entry:
\begin{align*}
\sum_k&\frac{(N-k-1)!k!(N+i-1)!}{i!(N-i-1)!(N+k-1)!(i-k)!}
\frac{(-1)^{j+N-1}(N+j-1)!}{j!(j-k)!(N-j-1)!k!}\\
&=\frac{(N+i-1)!(N+j-1)!}{i!(N-i-1)!j!(N-j-1)!}\sum_k\frac{(N-k-1)!}{(N+k-1)!(i-k)!}
\frac{(-1)^{j+N-1}}{(j-k)!}\\
&=\binom{N-1}{i}\binom{N+j-1}{j}\sum_{0\le k\le \min\{i,j\}}
(-1)^{j+N-1}\binom{N+i-1}{i-k}\binom{N-k-1}{j-k}.
\end{align*}
From this representation, it is clear that this is an integer. This was a question which was addressed in the affirmative in \cite{Richardson14}.

\section{$q$-analogues}

In this section we present $q$-analogues. Define $(q)_n:=(1-q)(1-q^2)\dots(1-q^n)$, and
\begin{equation*}
\gauss nk:=\frac{(q)_{n}}{(q)_{k}(q)_{n-k}};
\end{equation*}
these definitions are standard, see \cite{GaRa04}. Then we have the following results.

\begin{equation*}
L_{i,j}=\frac{(q)_i(q)_i(q)_{2j}}{(q)_{i+j}(q)_{i-j}(q)_{j}(q)_{j}},
\end{equation*}

\begin{equation*}
U_{i,j}=\frac{(-1)^iq^{i(3i-1)/2}(1+q^{i})(q)_{j}(q)_{j}(q)_{i}(q)_{i}}{(q)_{i+j}(q)_{j-i}(q)_{2i}}
\quad\text{for}\ i\ge1,\quad U_{0,j}=1,
\end{equation*}

\begin{equation*}
L^{-1}_{i,j}=\frac{q^{i(i-1)/2}(-1)^{i-j}(q)_{i}(q)_{i}(q)_{i+j-1}}{(q)_{2i-1}(q)_{i-j}}
\quad\text{for}\ j<i,\quad L^{-1}_{i,i}=1,
\end{equation*}

\begin{equation*}
U^{-1}_{i,j}=\frac{(-1)^iq^{-j^2-ji+i(i+1)/2}(q)_{j+i-1}(q)_{2j} (q)_{i}(q)_{i}}{(q)_{j-i}(q)_{j}(q)_{j-1}}
\quad\text{for}\ j>i,
\end{equation*}
\begin{equation*}
U^{-1}_{i,i}=\frac{(-1)^i q^{i(3i+1)/2} (q)_{2i}(q)_{2i}}{(q)_{i}(q)_{i}(q)_{i}(q)_{i}(1+q^{i})}
\quad\text{for}\ i\ge1,\quad U^{-1}_{0,0}=1,
\end{equation*}

\begin{equation*}
A_{i,j}=\frac{(-1)^{i-j}q^{(i+j+3)(i-j)/2+N(j-i)}(q)_{N-j-1}(q)_{j}(q)_{N+i-1}}{(q)_{N-i-1}(q)_i(q)_{N+j-1}(q)_{i-j}},
\end{equation*}

\begin{equation*}
B_{i,j}=\frac{(-1)^{j+N-1}q^{i^2+j(j+3)/2-Nj-N(N-1)/2}(q)_{N+j-1}}{(q)_j(q)_{j-i}(q)_{N-j-1}(q)_i},
\end{equation*}

\begin{equation*}
A_{i,j}^{-1}=\frac{q^{(i-j)(i-N+1)}(q)_{N-j-1}(q)_{N+i-1}(q)_j}{(q)_{N-i-1}(q)_{N+j-1}(q)_i(q)_{i-j}},
\end{equation*}

\begin{equation*}
B_{i,j}^{-1}=\frac{(-1)^{j+N-1}q^{j(j+1)/2-(N-j-1)i-N(N-1)/2}(q)_{j-i}(q)_{N+i-1}}{(q)_{N-i-1}(q)_j(q)_i}.
\end{equation*}

Note that for $q\to1$, we get the previous formul\ae. We do not discuss proofs here, since
Zeilberger's algorithm (aka WZ-theory)~\cite{PeWiZe96} proves all these results (which were obtained by guessing),
using a computer algebra system (such as, e.~g., Maple).

\medskip

\textsc{Remark.} Richardson's decomposition $S=GMG$ still holds when all binomial coefficients are replaced by the corresponding Gaussian $q$-binomial coefficients.

\bibliographystyle{plain}

\begin{thebibliography}{1}
	
	\bibitem{GaRa04}
	G.~Gasper and M.~Rahman.
	\newblock {\em Basic Hypergeometric Series}.
	\newblock Cambridge University Press, Cambridge, 2004.
	
	\bibitem{Gessel92}
	I.~M. Gessel.
	\newblock Super ballot numbers.
	\newblock {\em J. Symbolic Computation}, 14:179--194, 1992.
	
	\bibitem{GeLaFr03}
	I.~M. Gessel, P.~J. Larcombe, and D.~R. French.
	\newblock On the identity of von {Szily}: original derivation and a new proof.
	\newblock {\em Util. Math.}, 64:167--181, 2003.
	
	\bibitem{HeLiMa14}
	S.~Heubach, N.~Y. Li, and T.~Mansour.
	\newblock A garden of $k$-{Catalan} structures.
	\newblock {\em preprint}, 2014.
	
	\bibitem{PeWiZe96}
	M.~Petkov\v{s}ek, H.~Wilf, and D.~Zeilberger.
	\newblock {\em $A=B$}.
	\newblock A.K. Peters, Ltd., 1996.
	
	\bibitem{Richardson14}
	T.~M. Richardson.
	\newblock The reciprocal {P}ascal matrix.
	\newblock math.CO:arXiv:1405.6315, 2014.
	
	\bibitem{Szily94}
	K.~von Szily.
	\newblock {\"Uber die Quadratsummen der Binomialcoefficienten}.
	\newblock {\em Ungar. Ber.}, 12:84--91, 1894.
	
\end{thebibliography}

\end{document}